\documentclass[11pt]{article}

\title{  $\psi$-Poisson, $q$-Cigler ,$\psi$-Dobinski, $\psi$-Rota and $\psi$-coherent states-
with Cigler`s Remark on simplicity}
\author{A.K.Kwa\'sniewski\\  
\\ Higher School of Mathematics and Applied Informatics\\
PL - 15-021 Bialystok , ul.Kamienna 17,  Poland
\\e-mail: kwandr@uwb.edu.pl}

\usepackage{amsmath,amsthm}

\chardef\bslash=`\\ 
\begin{document}
\maketitle\
\begin{abstract}
The Cigler simple derivation of the $q$-Carlitz- Dobinski formula
is recalled and it is noticed that the formula  may be interpreted
as the average of powers of random variable $X_q$ with the
$q$-Poisson distribution. In parallel new $q$-Cigler-Dobinski and
$psi$-Carlitz-Dobinski formulas are introduced .
\end{abstract}

At first let us anticipate with $\psi$- remark. $\psi$ denotes an
extension of  $\langle\frac{1}{n!}\rangle_{n\geq 0}$ sequence to
quite arbitrary one ("admissible") and the specific choices are
for example : Fibonomialy  -extended ($\langle F_n \rangle$ -
Fibonacci sequence ) $\langle\frac{1}{F_n!}\rangle_{n\geq 0}$   or
just "the usual" $\langle\frac{1}{n!}\rangle_{n\geq 0}$ or Gauss
$q$-extended $\langle\frac{1}{n_q!}\rangle_{n\geq 0}$ admissible
sequences of extended umbral operator calculus - see more below.
With such an extension we may $\psi$-mnemonic repeat with exactly
the same simplicity and beauty what was done by Rota forty years
ago.
 Forty years ago Gian-Carlo Rota  \cite{1}  proved the exponential
generating function for Bell numbers $B_n$ to be of the form

\begin{equation}\label{eq1}
     \sum_{n=0}^\infty \frac {x^n}{n!}(B_n)= \exp(e^x-1)
\end{equation}
 using the linear functional  \textit{L } such that
\begin{equation}\label{eq2}
     L(X^{\underline{n}})=1  , \qquad  n\geq 0
\end{equation}\label{eq3}
Then Bell numbers (see: formula (4)  in \cite{1})  are defined by
\begin{equation}\label{eq4}
          L(X^n)=B_n  ,\qquad   n \geq 0
\end{equation}
The above formula is exactly the Dobinski formula \cite{2} if $L$
is interpreted as the average functional for the random variable
$X$  with the Poisson distribution with $L(X) = 1$. It is Blissard
calculus inspired umbral formula \cite{1} .
    Recently an interest to Stirling numbers and consequently
to Bell numbers was revived among "$q$-coherent states physicists"
\cite{3,4,5}. Namely the expectation value with respect to
coherent state $|\gamma >$ with $|\gamma| = 1$ of the $n$-th power
of the number of quanta operator is "just" the $n$-th Bell number
$B_n$ and the explicit formula for this expectation number of
quanta is "just" Dobinski formula \cite{3},(*). The same is with
the $q$-coherent states case \cite{3} i.e. the expectation value
with respect to $q$-coherent state $|\gamma> $ with  $|\gamma| =
1$ of the $n$-th power of the number operator is the $n$-th
$q$-Bell number \cite{6,3} defined as the sum  of $q$-Stirling
numbers  $ \Big\{{n \atop k}\Big\}_q $  due to Carlitz as in
\cite{6,3,4,5}. Note there then that  for  standard Gauss
$q$-extension $ x_q $ of number  $ x $ we have
\begin{equation}\label{eq5}
x_q^n=\sum_{k=0}^{n}\Big\{{n \atop k}\Big\}_q  x_q^{\underline k}
\end{equation}
Hence the expectation value with respect to $q$-coherent state $
|\gamma> $ with  $ |\gamma|=1 $ of the  $n$-th power of the number
operator  is exactly the popular $q$-Dobinski formula . It can be
given via (3) Blissard calculus inspired umbral formula form and
may be treated as definition of $ B_n(q)$
\begin{equation}\label{eq6}
         L_q(X_q^n)=B_n(q)  ,\qquad n\geq 0 .
\end{equation}
due to the fact that  linear functional  $L_q$    interpreted as
the average functional for the random variable  $X_q$  with the
$q$-Poisson distribution with  $L_q(X_q )= 1$  satisfies
\begin{equation}\label{eq7}
L_q(X_q^{\underline{n}})=1  , \qquad  n\geq 0.
\end{equation}
We arrive to this simple conclusion using Jackson derivative
difference operator in place of  $D = d/dx$  in $q$ =1 case and
the power series generating function $G(t)$  for $q$-Poisson
probability distribution:
\begin{equation}\label{eq8}
p_k=[exp_q\lambda]^{-1}\frac{\lambda^k}{k_q!},
G(t)=\sum_{n\geq0}p_k t^k ,
\end{equation}
\begin{equation}\label{eq9}
p_n=[\frac{\partial_q^n G(t)}{n_q!}]_{t=0} , [\partial_q
G(t)]_{t=1} = 1  for   \lambda = 1.
\end{equation}
There are many $q$-extensions of Stirling numbers according to
their weighted counting interpretation. For example $ w(\pi) =
q^{cross(\pi)} , w(\pi) = q^{inv(\pi)}$  from \cite{7} gives after
 being summed over  the  set of  $k$-block partitions  the
 Carlitz $q$-Stirling numbers  or $w(\pi) = q^{nin(\pi)} $
from \cite{8} gives rise to Carlitz-Gould $q$- Stirling numbers
after being summed over  the  set of  $k$-block partitions
or with  $w(\pi) = q^{i(\pi)}$ in \cite{9} - we arrive at another
combinatorial interpretation of $q$-extended Stirling numbers.
$q$-Stirling numbers much different from Carlitz $q$-ones were
introduced in the reference \cite{10} from where one infers \cite{11}
the \textit{cigl}-analog of (5) .  Let $\Pi$ denotes  the lattice of all
partitions of the set $\{0,1,..,n-1\}$.Let  $\pi\in \Pi $  be
represented by blocks   $ \pi =\{B_o ,B_1 ,...B_i ,...\}$ , where
$B_o$ is the block containing zero:  $0\in B_o$. The weight
adapted by Cigler  defines  weighted partitions` counting
according to the content of $B_o$. Namely $w(\pi)= q^{cigl(\pi)},
cigl(\pi)=\sum_{l\in {B_0}}l $,
$\sum_{\pi\in A_{n,k}}{q^{cigl(\pi)}}\equiv\Big\{{n\atop k}\Big\}_q $
therefore $ \sum_{\pi\in \Pi}{q^{cigl(\pi)}}\equiv{B_n(q)}$.
Here $A_{n,k}$ stays for subfamily of all $k$-block partitions.
With the above relations one has defined the \textit{cigl}-$q$-Stirling
 and the \textit{cigl}-$q$-Bell numbers.
The \textit{cigl}-$q$-Stirling numbers of the second kind  are
expressed in terms of $q$-binomial coefficients and $q =1$
Stirling numbers of the second kind \cite{10} . These are new
$q$-Stirling numbers. The corresponding \textit{cigl}-$q$-Bell
numbers recently have been equivalently defined via
\textit{cigl}-$q$-Dobinski formula \cite{11} $L(X_q^n)=B_n(q),
\qquad n\geq 0 , \qquad X_q^n\equiv X(X+q-1)...(X-1+q^{n-1}) $
interpreted as  the average of this specific  $n-th$
\textit{cigl}-$q$-power random variable $ X_q^{n} $ with the $q =
1$ Poisson distribution such that $ L(X)=1 $. To this end note
that in \cite{12} , \cite{13} a family of  the so called
$\psi$-Poisson processes was introduced. The corresponding choice
of the function sequence $\psi$ leads  to the $q$-Poisson process.
Accordingly the extension of Dobinski formula with its elementary
essential content and context to general case of $\psi$- umbral
instead $q$-umbral calculi case only - is automatic in view of  an
experience from \cite{12} , \cite{13} (see corresponding earlier
references there and necessary definitions). At first what you do
is to replace index $q$ by $\psi$ in formulas (3), (4),...,(8) .
Then you have got started problems with not easy combinatorial
interpretation if at all and... etc. $\psi$-Stirling numbers and
$\psi$-Bell numbers are being then defined by (4) and (3)
correspondingly with $q$ replaced by $\psi$. We get used to write
these extensions in mnemonic convenient
 upside down notation  \cite{12} , \cite{13}
\begin{equation}\label{eq10}
\psi_n\equiv n_\psi , x_{\psi}\equiv \psi(x)\equiv\psi_x ,
 n_\psi!=n_\psi(n-1)_\psi!, n>0 ,
\end{equation}
\begin{equation}\label{eq11}
x_{\psi}^{\underline{k}}=x_{\psi}(x-1)_\psi(x-2)_{\psi}...(x-k+1)_{\psi}
\end{equation}
\begin{equation}\label{eq12}
x_{\psi}(x-1)_{\psi}...(x-k+1)_{\psi}= \psi(x)
\psi(x-1)...\psi(x-k-1) .
\end{equation}
You may consult for further development and use of this notation
\cite{12} , \cite{13} and references therein.

(*) \textbf{Remark} based on the remark of Professor Cigler (
in private).\\
 The Katriel`s claim \cite{3} that his derivation of the Dobinski
 formula is the simplest possible may be confronted with the extremely
 simple derivation by Cigler (see p.$104$ in \cite{14}). Note on the way that
 this derivation is ready for  $\psi$-extensions \cite{12,13} as it
as a matter of fact based on elementary properties of $GHW$
(Graves-Heisenberg-Weyl)algebra : see \cite{12,13} and references therein.\\
Namely , let $\hat{x}$ denotes the multiplication by $x$ operator
while $D$ denotes differentiation -both acting on the prehilbert
space $P$  of polynomials . Then from recursion for Stirling
numbers of the second kind and the identification (operators in
$P$):
$$ \hat{x}(D+1)\equiv \frac {1}{\exp{(x)}}(\hat{x}D)\exp{(x)}$$
one gets for exponential polynomials
$$\varphi_n(x)=\sum_{k=0}^{n}\Big\{ {n \atop k}\Big\}x^k$$
the formula which becomes Dobinski one for  $x=1$
$$ \varphi_n(x)= \frac {1}{\exp{(x)}}(\hat{x}D)^n \exp{(x)}$$
i.e.
 $$\varphi_n(x)= \frac {1}{\exp{(x)}}\sum_{0\leq k}\frac{k^n
 x^k}{k!}.$$

 The $q$-case as well as $\psi$-case is automatically retained with
 the mnemonics from \cite{12,13}.

As perfectly properly indicated to the present author by Professor
Cigler -the derivation of Dobinski or $q$-Dobinski formulas  does
not require an introduction of a completion of prehilbert space
$P$ and coherent states in order to derive the formulas. - Well.
Anyhow , the extremely simple and umbral-beautiful Cigler`s
derivation above \cite{14} is immediately represented and then
fruitfully interpreted via expectation values on broader grounds
of combinatorics applications.

\end{document}